\documentclass[a4paper,12pt,reqno,tbtags]{amsart}
\usepackage{amsmath,amssymb,amsthm}
\usepackage{epsfig,graphicx,verbatim,color}

\numberwithin{equation}{section}
\newtheorem{thm}[equation]{Theorem}

\newtheorem{qn}[equation]{Question}
\theoremstyle{definition}

\newcounter{mycount}

\newenvironment{numlist}{\begin{list}{\arabic{mycount}.}%
   {\usecounter{mycount}\labelwidth=1cm\itemsep 0pt}}{\end{list}}

\newcommand\urladdrx[1]{{\urladdr{\def~{{\tiny$\sim$}}#1}}}

\renewcommand\o{{\mathrm o}}

\def\TT{{\mathbb T}}

\def\sS{{\mathcal S}}

\def\sW{{\mathcal W}}

\def\qq{\qquad}
\def\q{\quad}
\def\eps{\epsilon}
\def\b{\beta}
\def\a{\alpha}

\def\l{\lambda}

\def\t{\theta}

\def\NN{{\mathbb N}}
\def\ZZ{{\mathbb Z}}
\def\RR{{\mathbb R}}

\def\Pr{{\mathbb P}}

\def\Zcp{\ZZ^2_{\mathrm{cp}}}

\def\ra{\rangle}

\def\om{\omega}
\def\De{\Delta}
\def\de{\delta}

\def\oo{\infty}

\def\pc{p_{\mathrm c}}
\def\pcvec{\vec p_{\mathrm c}}

\def\be{\begin{equation}}
\def\ee{\end{equation}}

\def\resp{respectively}

\def\ra{\rightarrow}

\newcommand\subs[1]{\subseteq_{#1}}

\begin{document}
\title[Three problems for the clairvoyant demon]
{Three problems\\for the clairvoyant demon}

\author{Geoffrey Grimmett}
\address{Statistical Laboratory, Centre for Mathematical Sciences,
Cambridge University, Wilberforce Road, Cambridge CB3 0WB, UK}
\email{g.r.grimmett@statslab.cam.ac.uk}
\urladdrx{http://www.statslab.cam.ac.uk/~grg/}
\date{25 March 2009, revised 28 June 2009}

\begin{abstract}
A number of tricky problems in probability are discussed, having in
common one or more infinite sequences of coin tosses, and a 
representation as a problem in dependent percolation. 
Three of these problems are
of `Winkler' type, that is, they ask about what can be achieved 
by a clairvoyant demon.  
\end{abstract}

\maketitle

\section{Introduction}\label{sec:intro}

Probability theory has emerged in recent decades
as a crossroads where many sub-disciplines of mathematical science
meet and interact. 
Of the many examples within mathematics, we mention (not in order):
analysis, partial differential equations, mathematical physics, measure theory, 
discrete mathematics, theoretical computer science, and number theory.
The International Mathematical Union and the Abel Memorial Fund
have recently accorded acclaim to probabilists. This process
of recognition by others has been too slow, and would have been slower without the efforts
of distinguished mathematicians including John Kingman.

JFCK's work looks towards both theory and applications.
To single out just two of his theorems: the subadditive ergodic theorem
\cite{King1,King1a} is a piece of mathematical perfection that has also 
proved rather useful in practice; his `coalescent' \cite{King2a,King2b} is
a beautiful piece of probability that is now a lynchpin of
mathematical genetics. John is also an inspiring and devoted lecturer,
who continued to lecture to undergraduates even as the Bristol Vice-Chancellor,
and the Director of the Isaac Newton Institute in Cambridge.
Indeed, the current author learned his measure and probability from
partial attendance at John's course in Oxford in 1970/71. 

To misquote Frank Spitzer \cite[Chap. 8]{Spitz}, we turn to a very down-to-earth
problem: consider an infinite sequence of light bulbs. The basic commodity of 
probability is an infinite sequence of coin tosses. Such a sequence has been studied
for so long, and yet there remain `simple to state' problems that
appear very hard. We present some of these problems here. 
Sections \ref{sec:cl-sch}--\ref{sec:words}
are devoted to three famous problems for the so-called clairvoyant demon,
a non-human being to whom is revealed the (infinite) 
realization of the sequence, and who
is permitted to plan accordingly for the future.

Each of these problems may be phrased as a geometrical problem
of percolation type. The difference with classical percolation \cite{G99} lies
in the \emph{dependence} of the site variables. Percolation is reviewed
briefly in Section \ref{sec:sitep}.  This article ends with two
short sections on related problems, namely: other forms
of dependent percolation, and the question of `percolation of words'.

\section{Site percolation}\label{sec:sitep}

We set the scene by reminding the reader of the classical `site percolation
model' of Broadbent and Hammersley \cite{BH57}. 
Consider a countably infinite, connected
graph $G=(V,E)$. To each `site' $v \in V$ we assign a Bernoulli
random variable $\om(v)$ with density $p$. That is, $\om=\{\om(v): v \in V\}$
is a family of independent, identically distributed
 random variables  taking the values $0$ and $1$ with
respective probabilities $1-p$ and $p$.
A vertex $v$ is called \emph{open} if $\om(v)=1$, and \emph{closed} otherwise.

Let $0$ be a given vertex, called the \emph{origin}, and let $\t(p)$
be the probability that the origin lies in an infinite open
self-avoiding path of $G$. It is clear that $\t$ is non-decreasing
in $p$, and $\t(0)=0$, $\t(1)=1$. The \emph{critical probability}
is given as
$$
\pc = \pc(G) := \sup\{p: \t(p)=0\}.
$$
It is a standard exercise to show that the value of $\pc$ does not
depend on the choice of origin, but only on the graph $G$.

One may instead associate the random variables 
with the \emph{edges} of the graph, rather than the
\emph{vertices}, in which case the process is termed `bond percolation'.
Percolation is recognised as a fundamental model for
a random medium. It is important in probability and statistical physics,
and it continues to be the source of beautiful
and apparently hard mathematical
problems, of which the most outstanding
is to prove that $\t(\pc)=0$ for the three-dimensional
lattice $\ZZ^3$.  Of the several recent accounts of the percolation
model, we mention
\cite{G99,G-pgs}. 

\begin{figure}[tbp]
\centering
\includegraphics[width=9cm]{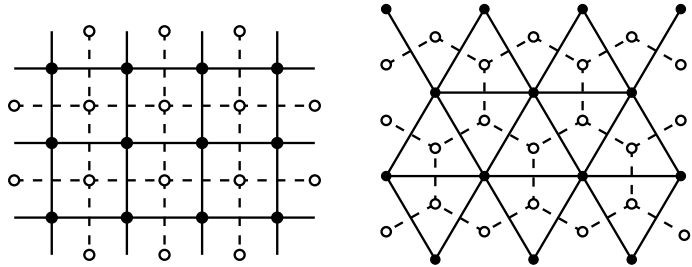}
\caption{The square lattice $\ZZ^2$ and the triangular lattice
$\TT$, with their dual lattices.}
\label{fig_lattices}
\end{figure}

Most attention has been paid to the case when $G$ is a crystalline 
lattice in two or more dimensions.  The current article is entirely concerned
with aspects of \emph{two-dimensional} percolation, particularly on the
square and triangular lattices illustrated in Figure \ref{fig_lattices}.
Site percolation on the triangular lattice has featured prominently
in the news in recent years, owing to the work of
Smirnov, Lawler--Schramm--Werner, and others on the relationship
of this model (with $p=\pc=\frac12$) to the process of random curves
in $\RR^2$ termed \emph{Schramm--L\"owner evolutions} (SLE),
and particularly the process denoted SLE$_6$.
See \cite{We07}.

When $G$ is a directed graph, one may ask 
about the existence of an infinite open \emph{directed} path
from the origin, in which case the process is referred to
as \emph{directed} (or \emph{oriented}) \emph{percolation}.

Variants of the percolation model are discussed in the following sections,
with the emphasis on models with site/bond variables that
are \emph{dependent}.

\section{Clairvoyant scheduling}\label{sec:cl-sch}

Let $G=(V,E)$ be a finite connected graph.
A symmetric random walk on $G$ is a sequence $X=(X_k: k=0,1,2,\dots)$
of
$V$-valued random variables satisfying, for all $v,w\in V$ and $k\ge 0$,
$$
P(X_{k+1}=w\mid X_k = v) = \begin{cases}
\dfrac 1 {\De_v} &\text{if } v \sim w,\\
0 &\text{if } v \nsim w,
\end{cases}
$$
where $\De_v$ is the degree of vertex $v$, and $\sim$
denotes the adjacency relation of $G$. Random walks on general graphs
have attracted much interest in recent years, see  
\cite[Chap.\ 1]{G-pgs} for example.

Let $X$ and $Y$ be independent random
walks on $G$ with distinct starting sites $x_0$, $y_0$, \resp.
We think of $X$ (\resp, $Y$) as describing the trajectory of
a particle labelled $X$ (\resp, $Y$) around $G$.
A clairvoyant demon is set the task of keeping the walks
apart from one another for all time. To this end, (s)he is 
permitted to schedule the walks in such a way that exactly one
walker moves at each epoch of time. Thus, the walks may be
delayed, but they are required to follow their prescribed trajectories.

More precisely, a \emph{schedule} is
defined as a sequence $Z=(Z_1,Z_2,\dots)$ in the 
space $\{X,Y\}^\NN$, and a given schedule
$Z$ is implemented in the following way. From the $X$ and $Y$ trajectories,
we construct the rescheduled walks $Z(X)$ and $Z(Y)$, where:

\begin{numlist}
\item If $Z_1 = X$, the $X$-particle takes one step
at time 1, and the $Y$-particle remains stationary. 
If $Z_1=Y$, it is the $Y$-particle that moves, and the $X$-particle 
that remains stationary. Thus,
\begin{alignat*}{3}
\text{if } Z_1=X\q &&\text{then}&\q Z(X)_1 = X_1,\ &Z(Y)_1=Y_0,\\
\text{if } Z_1=Y\q &&\text{then}&\q Z(X)_1 = X_0,\ &Z(Y)_1=Y_1.
\end{alignat*}

\item  Assume that, after time $k$, the $X$-particle has made
$r$ moves and the $Y$-particle $k-r$ moves, so that $Z(X)_k=X_r$
and $Z(Y)_k=Y_{k-r}$. 
\begin{alignat*}{3}
\text{If } Z_{k+1}=X\q &&\text{then}&\q 
  Z(X)_{k+1} = X_{r+1},\ &&Z(Y)_{k+1}=Y_{k-r},\\
\text{if } Z_{k+1}=Y\q &&\text{then}&\q 
   Z(X)_{k+1} = X_r,\ &&Z(Y)_{k+1}=Y_{k-r+1}.
\end{alignat*}

\end{numlist}

We call the schedule $Z$ \emph{good} if
$Z(X)_k\ne Z(Y)_k$ for all $k \ge 1$, and we say that the demon \emph{succeeds}
if there exists a good schedule $Z=Z(X,Y)$. (We overlook issues of measurability
here.)  The probability of success is
$$
\t(G) := \Pr(\text{there exists a good schedule}),
$$
and we ask: for which graphs $G$ is it the case that $\t(G)>0$?
This question was posed by Peter Winkler (see the discussion in 
\cite{CTW,Gacs2}). Note that the value of $\t(G)$ is independent
of the (distinct) starting points $x_0$, $y_0$.

The problem takes a simpler form when $G$ is the complete graph on some number, $M$ say,
of vertices. In order to simplify it still further, we add a loop to each
vertex. Write $V=\{1,2,\dots,M\}$, and $\t(M) := \t(G)$.
A random walk on $G$ is now a sequence of independent, identically
distributed points in $\{1,2,\dots,M\}$, each with the uniform distribution.
It is expected that $\t(M)$ is non-decreasing
in $M$, and it is clear by coupling that $\t(kM) \ge \t(M)$ for
$k \ge 1$. Also, it is not too hard to show that $\t(3)=0$.

\begin{qn}\label{winkler1}
Is it the case that $\t(M)>0$ for sufficiently large $M$?
Perhaps $\t(4)>0$?
\end{qn}

This problem has a geometrical formulation of percolation-type.
Consider the positive quadrant $\NN^2 = \{(i,j): i,j\ge 1,2,\dots\}$ of
the square lattice $\ZZ^2$. A \emph{path} is taken to be an infinite
sequence $(u_n,v_n)$, $n\ge 0$, with $(u_0,v_0)=(0,0)$
such that, for all $n\ge 0$,
$$
\text{either}\q (u_{n+1},v_{n+1}) = (u_n+1, v_n)
\q\text{or}\q 
(u_{n+1},v_{n+1}) =(u_n, v_n+1).
$$
With $X$, $Y$ the random walks
as above, we declare the vertex
$(i,j)$ to be \emph{open} if $X_i\ne Y_j$.  
It may be seen that the demon succeeds if and only if there exists
a path all of whose vertices are open.

Some discussion of this problem may be found in \cite{Gacs2}.
The law of the open vertices is $3$-wise independent but not $4$-wise independent
in the sense of Section \ref{sec:k-dep}.

The problem becomes significantly easier if paths are 
allowed to be undirected. For the totally undirected problem, it is
proved in \cite{BBS,W1} that there exists an infinite open path with
strictly positive probability if and only if $M \ge 4$.

\section{Clairvoyant compatibility}\label{sec:cl-comp}

Let $p \in(0,1)$, and let $X_1,X_2,\dots$ and $Y_1,Y_2\dots$ be
independent sequences of independent Bernoulli variables
with common parameter $p$. We say that a \emph{collision}
occurs at time $n$ if $X_n=Y_n=1$. The demon is now charged with
the removal of collisions, and to this end (s)he is
permitted to remove $0$s from the sequences. 

Let $\sW=\{0,1\}^\NN$, the set of singly-infinite sequences of $0$s and $1$s.
Each $w \in \sW$ is considered as a \emph{word} in an alphabet of two letters,
and we generally write $w_n$ for its $n$th letter.
For $w \in \sW$, there exists a sequence 
$i(w) = (i(w)_1,i(w)_2,\dots)$ of non-negative
integers such that $w = 0^{i_1}10^{i_2}1\cdots$, that is, there are
exactly $i_j=i(w)_j$ zeros between the $(j-1)$th and $j$th 
appearances of $1$.
For $x,y\in \sW$, we write $x \ra y$ 
if $i(x)_j \ge i(y)_j$ for $j\ge 1$. That is, $x \ra y$ if and only if $y$ may be obtained
from $x$ by the removal of 0s.

Two infinite words $v$, $w$ are said to be \emph{compatible}
if there exist $v'$ and $w'$ such that $v\ra  v'$,
$w \ra w'$ and $v'_nw'_n = 0$ for all $n$.
For given realizations $X$, $Y$, we say that the 
demon \emph{succeeds} if $X$ and $Y$ are compatible.
Write
$$
\psi(p) = \Pr_p(\text{$X$ and $Y$ are compatible}).
$$
Note that, by a coupling argument, $\psi$ is a non-increasing function.

\begin{qn}\label{winkler2}
For what $p$ is it the case that $\psi(p)>0$.
\end{qn}

It is easy to see as follows that $\psi(\frac12)=0$. When $p=\frac12$,
there exists a.s.\ an integer $N$ such that 
$$
\sum_{i=1}^N X_i > \tfrac12 N,\q \sum_{i=1}^N Y_i > \tfrac12 N.
$$
With $N$ chosen thus, it is not possible for the demon
to prevent a collision in the first $N$ values.
By working more carefully, one may see
that $\psi(\frac12-\eps)=0$ for small positive $\eps$;
see the discussion in \cite{gacs3}.

Peter G\'acs has proved in \cite{gacs3} that $\psi(10^{-400})>0$,
and he has noted that there is room for improvement.

\section{Clairvoyant embedding}\label{sec:words}
 
The clairvoyant demon's third problem stems from work on
long-range percolation of words (see Section \ref{sec:perc-words}). 
Let $X_1,X_2,\dots$ and $Y_1,Y_2,\dots$
be independent sequences of independent Bernoulli variables with
parameter $\frac12$. Let $M \in\{2,3,\dots\}$. The demon's task is
to find a monotonic embedding of the $X_i$ within the $Y_j$ in such a way
that the gaps between successive terms are no greater than $M$. 

Let $v,w \in \sW$. We say that $v$ is $M$-\emph{embeddable} in
$w$, and we write $v \subs M w$, if there exists an increasing sequence $(m_i: i \ge 1)$ 
of positive integers such that $v_i=w_{m_i}$ and $1\le m_i - m_{i-1}\le M$
for all $i \ge 1$. (We set $m_0=0$.) A similar definition is made for
\emph{finite} words $v$ lying in one of the spaces $\sW_n=\{0,1\}^n$, $n \ge 1$.

The demon succeeds in the above task if $X \subs M Y$,
and we let
$$
\rho(M) = \Pr(X \subs M Y).
$$
It is elementary that $\rho(M)$ is non-decreasing in $M$.

\begin{qn}\label{winkler3}
Is it the case that $\rho(M)>0$ for sufficiently large $M$?
\end{qn}

This question is introduced and discussed in \cite{GLR}, and partial but limited results proved.
One approach is to estimate the moments of the number $N_n(w)$
of $M$-embeddings of the finite word $w = w_1w_2\cdots w_n \in \sW_n$ 
within the random word $Y$.
It is elementary that $E(N_n(w)) = (M/2)^n$ for any such $w$,
and it may be shown that 
$$
\frac {E(N_n(X)^2)}{E(N_n(X))^2} \sim A_M c_M^n\qq \text{as } n\to\oo,
$$
 where
$A_M>0$ and $c_M > 1$ for $M \ge 2$.
That $E(N_n(w)) \equiv 1$ when $M=2$ is strongly suggestive that $\rho(2)=0$,
and this is part of the next theorem.

\begin{thm}\label{glr1}
\cite{GLR} 
We have that $\rho(2)=0$, and furthermore, for $M=2$,
\be
\Pr(w \subs 2 Y) \le \Pr(a_n \subs 2 Y)\qq\text{for all } w \in \sW_n,
\label{gen1}
\ee
where $a_n=0101\cdots$ is the \emph{alternating} word of length $n$.
\end{thm}

It is immediate that \eqref{gen1} implies $\rho(2)=0$ on noting that, for
any infinite
periodic word $\pi$, $\Pr(\pi \subs M Y) =0$ for all $M \ge 2$. One may 
estimate such probabilities more exactly through solving appropriate difference equations.
For example, $v_n(M) = \Pr(a_n \subs M Y)$ satisfies
\be
v_{n+1}(M) = (\a + (M-1)\b) v_n - \b(M-2\a) v_{n-1},\qq n \ge 1,
\label{gen2}
\ee
with boundary conditions $v_0(M)=1$, $v_1(M)=\a$. Here,
$$
\a  = 1-2^{-M},\q \b = 2^{-M}.
$$  
The characteristic polynomial associated with \eqref{gen2} is a quadratic with
one root in each of the disjoint intervals $(0,M\b)$ and $(\a,1)$. The larger
root equals $1-(1+\o(1))2^{1-2M}$ for large $M$, so that, in rough terms
$$
v_n(M) \approx (1-2^{1-2M})^n.
$$
Herein lies a health warning for simulators. One knows
that, almost surely, $a_n \not\subs M Y$ for large $n$, but one has to look
on scales of order $2^{2M-1}$ if one is to observe its extinction
with reasonable probability. 

One may ask about the `best' and `worst' words. Inequality \eqref{gen1} asserts that
an alternating word $a_n$ is the most easily embedded word when $M=2$. It is not known
which word is best when $M>2$. Were this a 
periodic word, it would follow that $\rho(M)=0$.
Unsurprisingly, the worst word is a constant word $c_n$ (of which there are of course two).
That is, for all $M \ge 2$,
$$
\Pr(w \subs M Y) \ge \Pr(c_n \subs M Y)\qq\text{for all } w \in \sW_n,
$$
where, for definiteness, we set $c_n=1^n\in \sW_n$.

Let $M=2$, so that the mean number $E(N_n(w))$ of embeddings of
any word of length $n$ is exactly 1 (as remarked above). A further argument is
required to deduce that $\rho(2)=0$. Peled \cite{RP} has made rigorous the following
alternative to that used in the proof of Theorem \ref{glr1}.
Assume that the word $w\in \sW_n$ satisfies $w \subs 2 Y$. For some small $c>0$,
one may identify (for most embeddings, with high probability) $cn$ 
positions at which the embedding may be altered, independently of each other. 
This gives $2^{cn}$ possible
`local variations' of the embedding.
It may be deduced that the probability of
embedding a word $w\in \sW_n$ is exponentially small in $n$, and also
$\rho(2)=0$.

The sequences $X$, $Y$ have been taken above with parameter $\frac12$. Little changes
in a more general setting. Let the two (respective) parameters be $p_X, p_Y\in (0,1)$. 
It turns out that the validity of the statement \lq\lq for all $M\ge 2$,
$\Pr(X \subs M Y) = 0$" is independent of the
values of $p_X$, $p_Y$. See \cite{GLR}.

A number of easier variations on Question \ref{winkler3} spring
immediately to mind, of which two are mentioned here. 
\begin{numlist}
\item Suppose the gap between the embeddings of $X_{i-1}$ and $X_{i}$
must be bounded above by some $M_i$. How slow a growth on the $M_i$ suffices that
the embedding probability be strictly positive? [An elementary bound
follows by the Borel--Cantelli lemma.]
\item Suppose that the
demon is allowed to look only boundedly into the future. How much clairvoyance
may (s)he be allowed without the embedding probability becoming
strictly positive?
\end{numlist}
Further questions (and variations thereof)
have been proposed by others.
\begin{numlist}
\item
In a `penalised embedding' problem, we are permitted mismatches by
paying a (multiplicative)  penalty $b$ for each. What is the cost of the `cheapest' penalised embedding of the first
$n$ terms, and what can be said as $b \to\oo$? [Erwin Bolthausen]
\item What can be said if we are required to embed only the 1s? That is, 
a `1' must be matched to a `1', but a `0' may be matched to either
`0' or `1'. [Simon Griffiths]
\item The above problems may be described as embedding $\ZZ$ in $\ZZ$.
In this language, might it be possible to embed $\ZZ^m$ in $\ZZ^n$ for some $m,n\ge 2$?
[Ron Peled]
\end{numlist}
 
Question \ref{winkler3} may be expressed as a geometrical problem of percolation type.
With $X$ and $Y$ as above, we declare the vertex $(i,j)\in\NN^2$ to be \emph{open} if
$X_i=Y_j$.  A \emph{path} in $\NN^2$ is defined  as an
infinite sequence $(u_n,v_n)$, $n\ge 0$,
of vertices such that: 
$$
(u_0,v_0)=(0,0),\q (u_{n+1},v_{n+1})=(u_n+1, v_n+d_n),
$$
for some $d_n$ satisfying $1\le d_n\le M$. 
It is easily seen that $X \subs M Y$ if and only if
there exists a path all of whose vertices are open.
(We declare $(0,0)$ to be open.)

\begin{figure}[tbp]
\centering
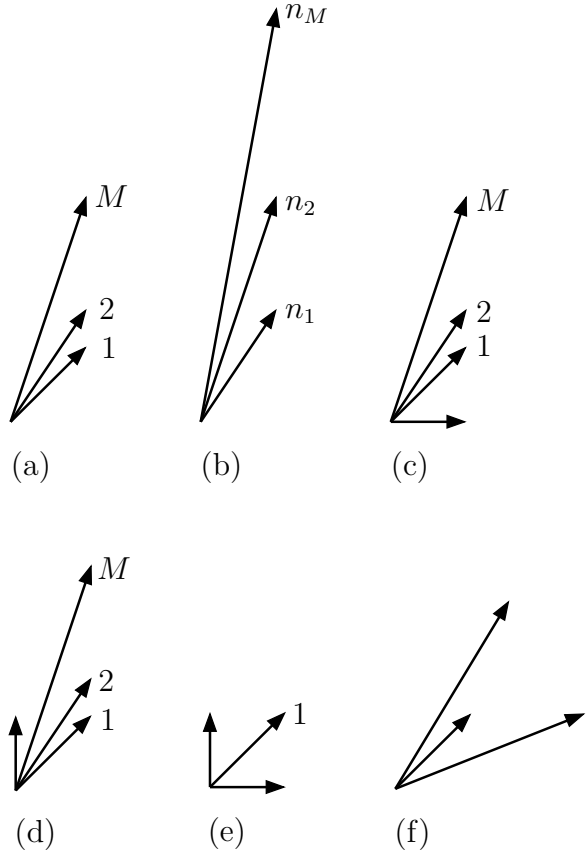
\caption{Icons describing a variety of embedding problems.}
\label{fig-icons}
\end{figure}

With this formulation in mind, the 
above problem may be represented by the icon 
at the top left of Figure \ref{fig-icons}. The further icons of that figure
represent examples of problems of similar type. 
Nothing seems to be known about these 
except:
\begin{numlist}

\item the argument of Peled \cite{RP} may be applied to problem (b)
with $M=2$ to obtain that $\Pr(w\subs 2 Y)=0$ for all $w \in \sW$, 

\item problem (e) is easily seen to be trivial. 

\end{numlist}

It is, as one might expect, much easier 
to embed words in two dimensions than in
one, and indeed this may be done along a path of $\ZZ^2$
that is directed in the north--easterly direction.
This statement may be made more 
precise as follows. Let $Y=(Y_{i,j}: i,j\ge 1)$ be a two-dimensional
array of independent Bernoulli variables 
with parameter $p \in (0,1)$, say. A word $v \in \sW$
is said to be $M$-embeddable in $Y$, written $v \subs M Y$,
if there exist strictly increasing sequences $(m_i: i \ge 1)$,
$(n_i: i \ge 1)$ of positive integers 
such that $v_i=Y_{m_i,n_i}$ and $1\le (m_i - m_{i-1})
+(n_i-n_{i-1})\le M$
for all $i \ge 1$. (We set $m_0=n_0=0$.)  
The following answers a question of \cite{deL}.
\emph{Note added at revision}: A related result has been discovered
independently in \cite{BSS}.

\begin{thm}\label{grg2}
\cite{G-pgs} Suppose $R \ge 1$ is such that $1 - p^{R^2} - (1-p)^{R^2}> \pcvec$, the critical probability
of directed site percolation on $\ZZ^2$. With strictly positive probability, 
every infinite word $w$ satisfies $w \subs {5R} Y$.
\end{thm}

The identification of the set of words that are $1$-embeddable in the two-dimensional
array $Y$, with positive probability, is much harder. This is a
problem of \emph{percolation of words}, and the results to
date are summarised in Section \ref{sec:perc-words}. 

\begin{proof}
We use a block argument. Let $R \in \{2,3,\dots\}$. For
$(i,j) \in \NN^2$, define the block $B_R(i,j) = ((i-1)R, iR]
\times ((j-1)R, jR] \subseteq \NN^2$. On the graph of blocks,
we define the (directed)  relation $B_R(i,j)\to B_R(m,n)$
if $(m,n)$ is either $(i+1,j+1)$ or $(i+1,j+2)$. By drawing a picture,
one sees that the ensuing directed graph is isomorphic to $\NN^2$
directed north--easterly. Note that the $L^1$-distance between two vertices
lying in adjacent blocks is no more than $5R$.

We call a block $B_R$ \emph{good} if it contains at least one 0 and at
least one 1. It is trivial that
$$
\Pr_p(B_R \text{ is good}) = 1 - p^{R^2} - (1-p)^{R^2}.
$$
If the right side exceeds the critical probability of directed site
percolation on $\ZZ^2$, then there is a strictly positive probability of
an infinite directed path of good blocks 
in the block graph, beginning at $B_R(1,1)$. Such a path contains
$5R$-embeddings of all words.
\end{proof}

The problem of clairvoyant embedding is connected to a question concerning
isometries of random metric spaces discussed in \cite{RP2}. In broad terms,
two metric spaces $(S_i,\mu_i)$, $i=1,2$, 
are said to be `quasi-isometric' (or `rough-isometric')
if their metric structure is the same up to multiplicative and additive
constants. That is, there exists a mapping $T: S_1 \to S_2$ and positive constants
$M$, $D$, $R$ such that:
$$
\frac 1M \mu_1(x,y)-D \le \mu_2(T(x),T(y)) \le M\mu_1(x,y)+D,
\qq x,y\in S_1,
$$
and, for $x_2\in S_2$, there exists $x_1\in S_1$ with
$\mu_2(x_2,T(x_1))\le R$.

It has been asked whether two Poisson process on the line, viewed as random
sets with metric inherited from $\RR$, are quasi-isometric. This question
is open at the time of writing.  
A number of related results are proved in \cite{RP2},
where a history of the problem may be found also. 
It turns out that the above question
is equivalent to the following. 
Let $X=(\dots,X_{-1},X_0,X_1,\dots)$ be a sequence of independent
Bernoulli variables with common parameter $p_X$. The 
sequence $X$ generates a random metric
space with points $\{i:X_i=1\}$ and metric inherited from $\ZZ$.
Is it the case that two independent sequences $X$ and $Y$ generate quasi-isometric
metric spaces? 
A possibly important difference between this problem and clairvoyant embedding is
that quasi-isometries of metric subspaces of $\ZZ$ need not be monotone.

\section{Dependent percolation}\label{sec:k-dep}

Whereas there is only one type of independence, there are many types
of dependence, too many to be summarised here. We mention just three
further types of dependent percolation in this section, of which the first (at least)
arises in the context of processes in random environments. In each, the
dependence has infinite range, and in this sense these problems 
have something in common with those treated in Sections 
\ref{sec:cl-sch}--\ref{sec:words}. 

For our first example, let $X = \{X_i: i \in \ZZ\}$ be independent, identically
distributed random variables taking values in $[0,1]$. Conditional
on $X$, the vertex $(i,j)$ of $\ZZ^2$ is
declared \emph{open} with probability $X_i$, and different vertices
receive (conditionally) independent states. 
The ensuing measure possesses a dependence that
extends without limit in the vertical direction. If the law $\mu$ of $X_0$
places probability both below and above $\pc$, there exist (almost surely)
vertically-unbounded domains that consider themselves subcritical, and others
that consider themselves supercritical. It depends on
the choice of $\mu$ whether or not the process possesses infinite
open paths, and necessary and sufficient conditions have proved elusive. 
The most successful technique for dealing with such problems seems to
be the so-called `multiscale analysis'. This leads to sufficient conditions under
which the process is subcritical (\resp, supercritical). 
See \cite{Klein1,Klein2}.

There is a variety of models of physics and applied probability
for which the natural random environment
is exactly of the above type. Consider, for example, the contact model in $d$ dimensions
with recovery rates $\de_x$ and infection rates $\l_e$, see \cite{Lig1,Lig2}.
Suppose that the environment is randomised
through the assumption that the $\de_x$ (\resp, $\l_e$) are independent and identically distributed.
The graphical representation of this model may be viewed as a
`vertically directed' percolation model on $\ZZ^d \times [0,\oo)$, in which the intensities
of infections and recoveries are dependent in the vertical direction.
See \cite{Andjel,BDS91,NewmanV} for example.

Vertical dependence arises naturally in certain models of statistical
physics also, of which we present one example. The `quantum Ising model'
on a graph $G$ may be formulated as a problem in stochastic
geometry on a product space of the form $G \times [0,\b]$, where $\b$ is
the inverse temperature. A fair bit of work has been done on
the quantum model in a random environment, that is, when its
parameters vary randomly around different vertices/edges of $G$.
The corresponding stochastic model on $G \times [0,\b]$ 
has `vertical dependence' of infinite range.
See \cite{BjG,GOS}. 

It is easy to adapt the above structure to provide dependencies in
both horizontal and vertical directions, although the
ensuing problems may be considered (so far) to have greater mathematical than physical interest. 
For example, consider bond percolation
on $\ZZ^2$, in which the states of horizontal edges are correlated thus,
and similarly those of vertical edges. A related
three-dimensional system has
been studied by Jonasson, Mossel and Peres \cite{JMP}.
Draw planes in $\RR^3$ orthogonal to the $x$-axis, 
such that they intersect the $x$-axis at points
of a Poisson process with given intensity $\lambda$.
Similarly, draw independent families of
planes orthogonal to the $y$- and $z$-axes. These three
families define a `stretched' copy of $\ZZ^3$. An edge of this stretched
lattice, of length $l$, is declared to be open with probability $e^{-l}$,
independently of the states of other edges. It is proved in \cite{JMP}
that, for sufficiently large $\lambda$, there exists (a.s.)
an infinite open directed percolation cluster that is
transient for simple random walk. The method of proof is interesting,
proceeding as it does by the method of `exponential intersection tails'
(EIT) of \cite{BPP}.
When combined with an earlier argument of H\"aggstr\"om, this
proves the existence of a percolation phase transition for the model.

The method of EIT is invalid in two dimensions, because random
walk is recurrent on $\ZZ^2$. The corresponding percolation
question in two dimensions was answered using different means by 
Hoffman \cite{Hof}.

In our final example, the dependence comes without geometrical
information. Let $k \ge 2$, and call a family of random
variables \emph{$k$-wise independent} if any $k$-subset is independent.
Note that the vertex states in Section \ref{sec:cl-sch} are $3$-wise independent 
but not $4$-wise independent.

Benjamini, Gurel-Gurevich and Peled \cite{BGGP} have investigated various
properties of $k$-wise independent Bernoulli families, and in particular the
following percolation question. Consider the $n$-box $B_n=[1,n]^d$ in $\ZZ^d$ 
with $d \ge 2$, in which the measure governing the site variables
$\{\om(v): v \in B_n\}$
has local density $p$ and is $k$-wise independent. Let $L_n$ be the
event that two given opposite faces are connected by an open
path in the box. Thus, for large $n$, the probability
of $L_n$ under the product measure $\Pr_p$ has a sharp threshold
around $p=\pc(\ZZ^d)$. The problem is to find bounds on the
smallest value of $k$ such that the probability of $L$ is
close to its value $\Pr_p(L_n)$ under product measure. 

This question may be formalised as follows. Let $\Pi = \Pi(n,k,p)$
be the set of probability measures on $\{0,1\}^{B_n}$
that have density $p$ and are $k$-wise independent. Let
$$
\eps_n(p,k) = \max_{\Pr\in \Pi} \Pr(L_n) - \min_{\Pr\in\Pi} \Pr(L_n),
$$
and
$$
K_n(p) = \min\{k: \eps_n(p,k) \le \de\},
$$
where for definiteness we may take $\de=0.01$ as in \cite{BGGP}.
Thus, roughly speaking, 
$K_n(p)$ is a quantification of the amount of independence required 
in order that, for all $\Pr\in\Pi$,
$\Pr(L_n)$ differs from $\Pr_p(L_n)$
by at most $\de$. 

Benjamini, Gurel-Gurevich and Peled have proved, in an ongoing project, that
$K_n(p) \le c \log n$ when $d = 2$ and $p \ne \pc$ (and when $d>2$
and $p<\pc$), for
some constant $c=c(p,d)$.
They have in addition a lower bound for $K_n(p)$
that depends on $p$, $d$, and $n$, and goes to $\oo$ as $n\to\oo$.

\section{Percolation of words}\label{sec:perc-words}

Recall the set $\sW=\{0,1\}^\NN$ of words in the alphabet
comprising the two letters $0$, $1$.
Consider the site percolation process of Section \ref{sec:sitep} 
on a countably infinite connected graph $G=(V,E)$, and write $\om=\{\om(v): v\in V\}$ for
the ensuing configuration.
Let $v\in V$ and let $\sS_v$ be the set of all self-avoiding walks
starting at $v$. Each $\pi\in\sS_v$ is a path $v_0,v_1,v_2\dots$ with $v_0=v$.
With the path $\pi$ we associate the word $w(\pi)= \om(v_1)\om(v_2)\cdots$,
and we write $\sW_v=\{w(\pi): \pi\in \sS_v\}$ for the set of
words `visible from $v$'. 
The central question of site percolation concerns the probability
that $\sW_v \ni 1^\oo$, where $1^\oo$ denotes the infinite word $111\cdots$. 
The so-called AB-percolation problem concerns the 
existence in $\sW_v$ of the infinite alternating word $01010\cdots$, see \cite{AW1}.

More generally, for given $p$, we ask which words lie in the 
random set $\sW_v$. 
Partial answers to this question may be found in three papers
\cite{BenjK,KSZ0,KSZ} of Kesten and co-authors
Benjamini, Sidoravicius, and Zhang, and they are summarised here as follows.

For $\ZZ^d$, with $p=\frac12$ and $d$ sufficiently large, we have from \cite{BenjK} that
$$
\Pr_{\frac12}(\sW_0=\sW) > 0,
$$
and indeed there exists (a.s.) some vertex $v$
for which $\sW_v=\sW$. Partial results are obtained for $\ZZ^d$
with edge-orientations in increasing coordinate directions.

For the triangular lattice $\TT$ and $p=\frac12$, we have from \cite{KSZ0} that
\be
\Pr_{\frac12}\left(\textstyle\bigcup_{v\in V} \sW_v 
\text{ contains almost every word}\right) = 1,
\label{gen3}
\ee
where the set of words seen includes all periodic words apart from $0^\oo$ and $1^\oo$.
The measure on $\sW$ can be taken in \eqref{gen3} as any 
non-trivial product measure. This extends the observation that AB-percolation
takes place at $p=\frac12$, whereas there is no infinite cluster
in the usual site percolation model. 

Finally, for the `close-packed' lattice $\Zcp$ obtained from $\ZZ^2$ by adding both
diagonals to each face,
$$
\Pr_p(\sW_0=\sW) > 0
$$
for $1-\pc < p < \pc$, with $\pc=\pc(\ZZ^2)$.
Moreover, every word is (a.s.) seen along some self-avoiding 
path in the lattice.
See \cite{KSZ}.

\section*{Acknowledgements}
The author acknowledges conversations with his co-authors Tom Liggett and Thomas Richthammer.
He profited from discussions with Alexander Holroyd while at the Department of Mathematics at the University
of British Columbia, and with Ron Peled and Vladas Sidoravicius while
visiting the Institut Henri Poincar\'e--Centre Emile Borel, 
both during 2008.
This article was written during a visit to the Section de Math\'ematiques
at the University of Geneva, supported by the Swiss National Science Foundation. 
The author thanks Ron Peled for his comments on a draft.

\def\PTRF{Probability Theory and Related Fields\ }
\bibliography{demon} \bibliographystyle{amsplain}

\end{document}